\documentclass[11pt,a4paper]{amsart}

\usepackage{amsmath}

\usepackage[dvipdfmx]{pict2e}
\usepackage[dvipdfmx]{hyperref}
\usepackage{amsfonts}

\usepackage{adjustbox}

\usepackage[abbrev]{amsrefs}

\makeatletter

\@addtoreset{equation}{section}
\makeatother

\DeclareMathOperator{\End}{End}
\DeclareMathOperator{\Int}{Int}
\DeclareMathOperator{\Ker}{Ker}

\DeclareMathOperator{\Ima}{Im}

\DeclareMathOperator{\vol}{vol}

\DeclareMathOperator{\Id}{Id}
\DeclareMathOperator{\prim}{prim}

\newcommand{\adjoint}{\dagger}
\newcommand{\deRham}{\mathrm{dR}}
\newcommand{\Rumin}{\mathrm{R}}
\newcommand\RN {\calE }

\newcommand{\bC}{\mathbb{C}}
\newcommand{\calE}{\mathcal{E}}

\newcommand{\BigWedge}%
{{\mathord{\adjustbox{valign=B,totalheight=.6\baselineskip}{$\bigwedge$}}}}

\theoremstyle{plain}
\newtheorem{defi}{Definition}[section] %section
\newtheorem{remark}[defi]{Remark}
\newtheorem{corollary}[defi]{Corollary}
\newtheorem{proposition}[defi]{Proposition}
\newtheorem{theorem}[defi]{Theorem}
\newtheorem*{theorem*}{Theorem}
\newtheorem{lemma}[defi]{Lemma}

\begin{document}
\title[Analytic Tonsions associated with the Rumin Complex]
{
Analytic Torsions associated with \\
the Rumin Complex
on Contact Spheres
}

\author{Akira Kitaoka}
\address{Graduate School of Mathematical Sciences \\ The University of Tokyo \\ 3-8-1 Komaba, Meguro, Tokyo 153-8914 Japan}
\email{akitaoka@ms.u-tokyo.ac.jp}

\begin{abstract}
We explicitly write down all eigenvalues of the Rumin Laplacian on the standard contact spheres,
and express the analytic torsion functions associated with the {Rumin complex}
in terms of the Riemann zeta function.  In particular, we find that the functions vanish at the origin and determine the analytic torsions.
\end{abstract}

\maketitle

\subsection*{Introduction}
Let $(M,H)$ be a compact contact manifold of dimension $2n+1$.
Rumin \cite{Rumin-94} introduced a complex $(\calE ^\bullet , d_{\Rumin }^\bullet)$,
which is
a resolution of the constant sheaf of $\mathbb{R}$
given by
a subquotient of the de Rham complex.
A specific feature of the complex is that
the operator $D = d_{\Rumin }^n \colon \calE ^{n} \to \calE ^{n+1}$ in `middle degree' is a second-order,
while
 $d_{\Rumin }^k \colon \calE ^{k} \to \calE ^{k+1} $ for $k \not = n$ are first order which are induced by
the exterior derivatives.
Let $a_k = 1 / \sqrt{|n-k|}$ for $k \not = n$ and $a_n = 1$.
Then, $(\calE ^\bullet , d_\RN ^\bullet)$,
where $d_\RN ^k = a_k d_{\Rumin}^k$,
is also a complex.
We call $(\calE ^\bullet , d_\RN ^\bullet)$ the {\em  Rumin complex.}
In virtue of the rescaling,  $ d_\RN ^\bullet$
satisfies K\"{a}hler-type identities on Sasakian manifolds \cite{Rumin-00}, which include the case of  spheres.

Let $\theta$ be a contact form of $H$
and
$J$ be an almost complex structure on $H$.  Then we may define a Riemann metric $g_\theta$ on $TM$ by extending the Levi metric $ d\theta (- ,J-)$ on $H$ (see \S1.1).  Following  \cite{Rumin-94}, we define the Rumin Laplacians $\Delta_\RN $ associated with $(\calE ^\bullet , d_\RN ^\bullet)$ and the metric $g_\theta$ by
\[
\Delta_\RN ^k
:=
\begin{cases}
( d_\RN d_\RN {}^{\adjoint })^2
+ (d_\RN {}^{\adjoint } d_\RN )^2 , & k \not = n, \, n+1 , \\
( d_\RN d_\RN {}^{\adjoint })^2
+ D^{\adjoint } D , & k=n, \\
D D^{\adjoint }
+ (d_\RN {}^{\adjoint } d_\RN )^2 , & k =n+1.
\end{cases}
\]
Rumin  showed that $\Delta_\RN $ have discrete eigenvalues with finite multiplicities.

In this paper, we determine
explicitly eigenvalues of $\Delta_\RN $
on the standard contact spheres $S^{2n+1}\subset \bC^{n+1}$.
The sphere $S^{2n+1}$ also admits an almost complex structure $J$ induced from the complex structure of $ \bC^{n+1}$.
We call $(S^{2n+1}, H, J)$ the {\em standard CR sphere} and simply denoted by $S^{2n+1}$.
To state our result we need to introduce notation for highest weight representations of the unitary group $U(n+1)$
which acts on $S^{2n+1}$.
The irreducible representations of  $U(n+1)$ are classified by the highest weights $\lambda=(\lambda_1,\lambda_2,\dots,\lambda_{n+1})$; the corresponding representation will be denoted by $V(\lambda)$.   Julg and Kasparov \cite{Julg-Kasparov-95} showed that
the complexification of $\calE ^k (S^{2n+1})$, as a $U(n+1)$-module, is decomposed into the irreducibles of the form
 \[
 \Psi_{(q,j,i,p)} := V (q, \underbrace{1, \ldots , 1 }_{j \text{ times }} , 0, \ldots , 0, \underbrace{-1 , \ldots , -1}_{i \text{ times}} , -p ).
 \]
See Proposition \ref{theo:Irreducible-decomposition-of-Rumin-complex} below for the relations between $k$ and $(q,i,j,p)$.
Since $\Delta_\RN $ commutes with the $U(n+1)$-action, it acts as a scalar on each $\Psi_{(q,j,i,p)}$.

%The above sentence has incomplete sentences.
\begin{theorem}
\label{theo:all eigenvalues of the Rumin Laplacian on sphere}
Let $S^{2n+1}$ be the standard CR sphere with the contact from $\theta = \sqrt{-1} (\bar{\partial}  - \partial) |z|^2$.Then, on the subspaces of the complexification of $\calE ^\bullet$
corresponding to the representations
$\Psi_{(q,j,i,p)}$,
the eigenvalue of $\Delta_\RN $ is
\[
\frac{\left( ( p +i )(q +n -i) + ( q +j )(p +n -j) \right)^2}{4 (n-i-j)^2 } .
\]
\end{theorem}

This theorem claims that the eigenvalues of $\Delta_\RN $ are determined by the highest weight.
This phenomenon also appears in the case of the Hodge-de Rham Laplacian $\Delta_{\deRham}$ on symmetric spaces $G / K$.
Ikeda and Taniguchi \cite{Ikeda-Taniguchi-78} showed that
on the subspaces of $k$-forms of $G/K$
corresponding to $V (\lambda )$,
the eigenvalue of $\Delta_{\deRham}$ is determined by $\lambda$.
It is a natural question to ask
whether
the eigenvalues of $\Delta_\RN $
on a contact homogeneous space $G/K$
are determined by the highest weight  of $G$.

Theorem \ref{theo:all eigenvalues of the Rumin Laplacian on sphere} unifies the follwoing  results on the eigenvalues of  Rumin Laplacians  on the spheres.
Julg and Kasparov \cite{Julg-Kasparov-95} determined the eigenvalues of $D^{\adjoint} D $.
Folland \cite{Folland-72} calculated the eigenvalue of the sub-Laplacian $\Delta_b$,
which agrees with $\Delta_\RN $ on $\calE ^{0}$.
Seshadri \cite{Seshadri-07} determined the eigenvalues of $d_\RN  d_\RN ^{\adjoint}$ on $\calE ^{1}$ in the case $S^3$.
{\O}rsted and Zhang \cite{Orsted-Zhang-05} determined eigenvalues of the Laplacian of the holomorphic and anti-holomorphic part of $d_{\Rumin}$ except for the ones containing $D$.

Note that {\O}rsted and Zhang  used $d_R$ in place of $d_\RN $. As a result, the eigenvalues of the Laplacian in their paper are not determined  by the highest weights.  This also explains the importance the scaling factor $a_k$.

We next introduce  the  analytic torsions and metrics of the Rumin complex $(\calE ^\bullet  , d_\RN ^\bullet)$ by following \cite{Bismut-Zhang-92 , Rumin-Seshadri-12}.
We define the {\em contact analytic torsion function} associated with $(\calE ^\bullet, d_\RN ^\bullet)$ by
\begin{equation}\label{kappa}
\kappa_\RN (s)  := \, \sum_{k=0}^{n} (-1)^{k+1} (n+1-k) \zeta (\Delta_\RN ^k ) (s),
\end{equation}
where $\zeta (\Delta_\RN ^k ) (s)$ is the spectral zeta function of $\Delta_\RN ^k$,
and the {\em contact analytic torsion} $ T_\RN $  by
\[
2\log T_\RN
=
\kappa_\RN ' (0) .
\]
Let $H^\bullet (\calE ^\bullet , d_\RN ^\bullet)$ be the cohomology of the Rumin complex.
We define the contact metric on $\det H^\bullet (\calE ^\bullet , d_\RN ^\bullet)$ by
\[
\| \quad \|_\RN
=
T_\RN
| \quad |_{L^2 (\calE ^\bullet) },
\]
where
the metric
$| \quad|_{L^2 (\calE ^\bullet) }$ is induced by $L^2$ metric on $\calE ^\bullet$
via identification of the cohomology classes by the harmonic forms on $\calE ^\bullet$.
%the above sentence is what you want to say.

Rumin and Seshadri \cite{Rumin-Seshadri-12} defined  another analytic torsion function $\kappa_{\Rumin}$ from
$(\calE ^\bullet , d_{\Rumin }^\bullet)$, which is different from $\kappa_\RN $ except in dimension 3.
In dimension 3, they showed that
$\kappa_{\Rumin} (0)$ is a contact invariant; so is $\kappa_\RN  (0)$.
Moreover, on 3-dimensional Sasakian manifolds with $S^1$-action, $\kappa_{\Rumin} (0) = 0 $ and
the contact analytic torsion and metric coincides with the Ray-Singer torsion $T_{\deRham }$ and the metric $\| \quad \|_{\deRham }$.

Our second main result is
\begin{theorem}
\label{theo:compute-precisely-kappa-on-the-sphere}
In the setting of Theorem \ref{theo:all eigenvalues of the Rumin Laplacian on sphere},
the contact analytic torsion function
on $S^{2n+1}$
is given by
\[
\kappa_\RN (s)
= \, - (n+1) (1+ 2^{2s+1} \zeta (2s ))
,
\]
where $\zeta$ is the Riemann zeta function.
In particular, we have
\begin{align}
& \kappa_\RN (0)
=
0
,
\label{eq:kappa-at-the-origin-on-the-sphere}
\\
& T_\RN
=
(4 \pi)^{n+1}
.
\label{eq:Rumin-analytic-torsion-on-the-sphere}
\end{align}
\end{theorem}

From (\ref{eq:kappa-at-the-origin-on-the-sphere}),
we see that
the metric $\| \quad\|_\RN $ on $S^{2n+1}$ is invariant under the constant rescaling $\theta \mapsto K \theta$.
The argument is exactly same as the one in \cite{Rumin-Seshadri-12}.

The fact that the representations determines the eigenvalues of $ \Delta_\RN$ causes several cancelations in the linear combination \eqref{kappa}, which greatly simplifies the computation of $\kappa_\RN (s)$.
We cannot get such a simple formula for the contact torsion function $\kappa_{\Rumin}$ of $(\calE ^\bullet, d_{\Rumin}^\bullet)$
for dimensions higher than 3.

Let us compare the contact analytic torsion with the Ray-Singer torsion on the standard spheres given by Weng and You \cite{Weng-You-96}.
Let $g_{\mathrm{std}}$ be the standard metric on $S^{2n+1}$.
They showed that the Ray-Singer torsion of $(S^{2n+1} , 4 g_{\mathrm{std}} )$ is $(4 \pi)^{n+1} /n!$.
The metric $4 g_{\mathrm{std}}$ agrees with the metric $g_\theta$ defined from the contact from $\theta = \sqrt{-1} (\bar{\partial}  - \partial) |z|^2$.
Since $(\calE ^\bullet  , d_\RN ^\bullet)$ and $ (\Omega^\bullet  , d )$ are resolutions of $\mathbb{R}$, there is a natural  isomorphism
$\det H^\bullet (\calE ^\bullet  , d_\RN ^\bullet)
\cong
\det H^\bullet (\Omega^\bullet  , d )
$, which turns out to be isometric for the $L^2$ metrics.
Therefore \eqref{eq:kappa-at-the-origin-on-the-sphere} gives
\begin{corollary}
\label{cor:Rumin-and-Ray-Singer-analytic-torsion-on-the-sphere}
In the setting of Theorem \ref{theo:all eigenvalues of the Rumin Laplacian on sphere},
we have
\[
T_\RN = n ! T_{\deRham }
\quad \text{and} \quad
\|\quad \|_\RN = n! \| \quad \|_{\deRham }.
\]
\end{corollary}

The paper is organized as follows.
In \S \ref{sec: The Rumin complex on the spheres},
we recall the definition and properties of the Rumin complex on $S^{2n+1}$,
and decompose $\calE ^{k}$ as a direct sum of the irreducible representation of $U(n+1)$.
In \S \ref{sec:The action of the Rumin derivative and the Reeb vector field},
we construct highest weight vectors,
and compute the actions of $d_{\Rumin }$ and the Lie derivative $\mathcal{L}_T$ with respect to the Reeb vector field $T$ on these vectors.
In \S \ref{subsec: L^2 norm of elements},
we calculate the $L^2$-norm of them.
Then, in \S \ref{subsec: Calculation of Rumin derivative and Laplacians},
we compute the eigenvalues of $\Delta_\RN $ for each irreducible component.
In \S \ref{sec:The analytic torsion function},
we calculate the contact analytic torsion function $\kappa_\RN $.

\subsection*{Acknowledgement}
The author is grateful to his supervisor Professor Kengo Hirachi for introducing this subject and for helpful comments.
He would also like to thank Yuya Takeuchi for carefully reading an earlier draft of this paper.
This work was supported by the program for Leading Graduate Schools,
MEXT, Japan.

\section{The Rumin complex}
\label{sec: The Rumin complex on the spheres}

\subsection{The Rumin complex on contact manifolds}

We call $(M,H)$ an orientable contact manifold of dimension $2n+1$
if $H$ is a subbundle of $TM$ of codimension $1$
and
there exists a 1-form $\theta$,
called a contact form,
such that
 $\Ker (\theta \colon TM \to \mathbb{R} ) = H $
and $\theta \wedge ( d \theta )^n \not = 0$.
The Reeb vector field of $\theta$ is the unique vector field $T$ satisfying $\theta (T) =1$ and $\Int_T d \theta =0$, where $\Int_T$ is the interior product with respect to $T$.

For $H$ and $\theta$,
we call $J \in \End (TM)$
an almost complex structure associated with $\theta$ if
$J^2 = - \Id $ on $H$, $JT=0$, and the Levi form $d \theta ( - , J - )$ is positive definite on $H$.
Given $\theta $ and $J$, we define a Riemannian metric $g$ on $TM$ by
\begin{align*}
g (X,Y) := d \theta (X , JY) + \theta (X) \theta (Y) \hspace{1cm} \text{ for } X,Y \in TM .
\end{align*}
Let $*$ be the Hodge star operator on $\wedge^\bullet T^*M$ with respect to $g$.

The Rumin complex \cite{Rumin-94} is defined on contact manifolds as follows.
We set $L := d \theta \wedge $
and $\Lambda := *^{-1} L *$,
which is the adjoint operator of $L$ with respect to the metric $g$ at each point.
We set
\begin{align*}
\BigWedge_{\prim}^{k} H^*&:=
\left\{
v \in \BigWedge^{k} H^*
\middle|
\Lambda v = 0
\right\},
\\
\BigWedge_{L}^{k} H^*
& :=
\left\{
v \in \BigWedge^{k} H^*
\middle|
L v = 0
\right\},
\\
\calE ^k & :=
\left\{
\begin{aligned}
& C^\infty \left( M,\BigWedge_{\prim}^{k} H^* \right) , \quad && k \leq n ,\\
& C^\infty \left( M, \theta \wedge\BigWedge_{L}^{k-1} H^* \right) , &&  k \geq n+1 .
\end{aligned}
\right.
\end{align*}
We embed $H^*$ into $T^*M$ as  the subbundle $ \{ \phi \in T^* M \mid \phi (T) = 0 \}$ so that
we can regard
\begin{align*}
\Omega_H^{k}  := & \, C^{\infty} \left( M, \BigWedge^{k} H^* \right)
\end{align*}
as a subspace of $\Omega^k$, the space of $k$-forms.
We define  $d_b \colon \Omega_H^{k}  \to \Omega_H^{k+1} $ by
\begin{align}
d_b \phi := & \, d \phi - \theta \wedge (\Int_T d \phi )
\label{eq: d_b-1}.
\end{align}
and then $D \colon \calE ^n \to \calE ^{n+1} $  by
\begin{align}
D = \theta \wedge ( \mathcal{L}_T + d_b L^{-1} d_b ),
\label{eq:D-0}
\end{align}
where $\mathcal{L}_T$ is the Lie derivative with respect to $T$,
and we use the fact that $L \colon \bigwedge^{n-1} H^* \to\bigwedge^{n+1} H^*$ is an isomorphism.

Let $P \colon \bigwedge^{k} H^* \to \bigwedge_{\prim}^{k} H^* $ be the fiberwise orthogonal projection with respect to $g$, which also defines a projection  $P \colon \Omega^k  \to \calE ^k  $.
We set
\[
d_{\Rumin }^k
:=
\begin{cases}
P \circ d &  \text{ on } \calE ^k,  \quad k \leq n-1, \\
D & \text{ on } \calE ^n  ,  \\
d  & \text{ on } \calE ^k, \quad  k \geq n+1.
\end{cases}
\]
Then $(\calE ^\bullet, d_{\Rumin }^\bullet)$ is a complex.
Let  $d_\RN ^{k} = a_k d_{\Rumin }^{k}$, where
$a_k = 1 / \sqrt{|n-k|}$ for $k \not = n$ and
$a_n = 1$.
We call $(\calE ^\bullet, d_\RN ^\bullet)$ the Rumin complex.

We define the $L^2$-inner product on $\Omega^{k}$ by
\[
(\phi , \psi ) := \int_M g (\phi , \psi ) d \vol_g
\]
and the $L^2$-norm on $\Omega^k $ by $\| \phi \| := \sqrt{(\phi , \phi )}$.
Since the Hodge star operator $*$ induces a bundle isomorphism
from $ \bigwedge_{\prim}^{k} H^*$ %\overset{\cong}{\to}
to $\theta \wedge \bigwedge_{L}^{2n-k} H^*$,
it  also induces a map $\calE ^{k}  \to \calE ^{2n+1-k} $.

Let $d_\RN ^{\adjoint }$ and $D^{\adjoint }$ denote the formal adjoint of $d_\RN $ and $D$, respectively for the $L^2$-inner product.
We define the forth-order Laplacians $\Delta_\RN  $ on $\calE ^k$ by
\[
\Delta_\RN ^k
 :=
\begin{cases}
( d_\RN ^{k-1} d_\RN ^{k-1} {}^{\adjoint })^2 + (d_\RN ^{k} {}^{\adjoint } d_\RN ^{k} )^2 , & k \not = n , n+1 , \\
\vspace{1mm}
( d_\RN ^{n-1} d_\RN ^{n-1} {}^{\adjoint })^2 + D^{\adjoint } D , & k=n, \\
D D^{\adjoint } + (d_\RN ^{n+1} {}^{\adjoint } d_\RN ^{n+1} )^2 , & k =n+1.
\end{cases}
\]
We call it the Rumin Laplacian \cite{Rumin-94}.
Since $*$ and $\Delta_\RN $ commute,
to determine the eigenvalue on $\calE ^\bullet  $,
it is enough to compute them on $\calE ^{k} $ for $k \leq n$.

\subsection{The Rumin complex on the CR spheres}

Let $S :=\{ z \in \bC^{n+1} \mid  |z|^2 =1 \}$
and
$\theta:= \sqrt{-1} (\bar{\partial } - \partial )|z|^2$.
(We will omit the dimension from $S^{2n+1}$  for the simplicity of the notation.)
The Reeb vector filed of $\theta$ is
\[
T = \frac{\sqrt{-1}}{2} \sum_{l=1}^{n+1}
\left(
z_l \frac{\partial }{\partial z_l } - \bar{z}_l \frac{\partial }{\partial \bar{z}_l }
\right).
\]
With respect to the standard almost complex structure $J$, we decompose the bundles defined in the previous subsection as follows:
\begin{align*}
H^{* 1,0} &:=  \{ v \in \bC H^*   \mid Jv = \sqrt{-1} v \} , \\
H^{* 0,1}& := \{ v \in \bC H^*  \mid Jv = - \sqrt{-1} v \}, \\
\BigWedge^{i, j} H^*& :=   \BigWedge^{i} H^{* 1,0} \otimes \BigWedge^{j} H^{* 0,1} , \\
\BigWedge_{\prim}^{i, j} H^*& :=   \left\{ \phi \in \BigWedge^{i, j} H^* \middle| \Lambda \phi =0 \right\} , \\
\Omega_H^{i, j}& := C^{\infty} \left( S, \BigWedge^{i, j} H^* \right) , \\
\calE ^{i, j}& := C^{\infty} \left( S, \BigWedge_{\prim}^{i, j} H^* \right).
\end{align*}
Then $d_b \Omega_H^{i , j} \subset \Omega_H^{i+1 , j} \oplus \Omega_H^{i , j+1} $.
We define $\partial_b \colon \Omega_H^{i , j}$
$\to$
$\Omega_H^{i+1 , j} $ and $\bar{\partial}_b \colon \Omega_H^{i , j} \to \Omega_H^{i , j+1 } $ by
\[
d_b = \partial_b + \bar{\partial}_b.
\]
Similarly, we decompose
\[
d_{\Rumin} = \partial_{\Rumin} + \bar{\partial}_{\Rumin}
,
\quad
d_\RN  = \partial_\RN  + \bar{\partial}_\RN
.
\]
In view of the Lefshetz primitive decomposition, we may rewrite \eqref{eq:D-0} as
\begin{equation}
D = \theta \wedge \left( \mathcal{L}_T - \sqrt{-1} (\partial_b + \bar{\partial}_b )
( \partial_b^{\adjoint } - \bar{\partial}_b^{\adjoint } )\right)
\label{eq:D-1}
\end{equation}
by using
$\partial_b^{\adjoint} = \sqrt{-1} [ \Lambda , \bar{\partial}_b]$ and
$\bar{\partial}_b^{\adjoint} = - \sqrt{-1} [ \Lambda , \partial_b]$.
Note that this equation holds on Sasakian manifolds.

We decompose $\calE ^{i,j} $ into a direct sum of irreducible representations of the unitary group $U(n+1)$.
Recall that
irreducible representations of $U (m)$ are parametrized by the highest weight $\lambda = (\lambda_1 , \ldots , \lambda_m ) \in \mathbb{Z}^m$ with $\lambda_1 \geq \lambda_2 \geq \cdots \geq \lambda_m$; the representation corresponding to $\lambda$ will be denoted by
$ V ( \lambda ) $.
To simplify the notation, we introduce the following notation:
for $a_1, \ldots $, $a_l \in \mathbb{Z}$ and $k_1, \ldots , k_l \in \mathbb{Z}$, $(\underline{a_1}_{k_1}, \cdots , \underline{a_l}_{k_l})$ denotes the $k_1 + \cdots + k_l $-tuple
whose first $k_1$ entries are $a_1$, whose next $k_2$ entries are $a_2$, etc.
For example,
\[
(\underline{1}_{3}, \underline{0}_{2}, \underline{-1}_{2} ) = (1, 1, 1, 0, 0, -1, -1 ).
\]
We note that $\underline{a}_{1}$ is $a$ and $\underline{a}_{0}$ is the zero tuple.

In \cite{Julg-Kasparov-95}, it is shown that
the multiplicity of $V (q, \underline{1}_{j}, \underline{0}_{n-1-i-j}, \underline{-1}_{i}, -p )$ in $\calE ^{s,t}$
is at most one.
Thus we may set
\begin{equation}
\Psi_{(q , j , i , p )}^{(s, t)}
:= \calE ^{s,t}
\cap V (q, \underline{1}_{j}, \underline{0}_{n-1-i-j}, \underline{-1}_{i}, -p )
.
\label{eq:irreducible_components_of_E}
\end{equation}

\begin{proposition}
{\rm
(\cite[Section 4(b)]{Julg-Kasparov-95})
}
\label{theo:Irreducible-decomposition-of-Rumin-complex}
The irreducible decomposition of  the $U(n+1)$-module $\calE ^{i,j} $ is given as follows:

\noindent
Case $\mathrm{I}:$
\[
\calE ^{0,0}
=
\bigoplus_{q \geq 0 , p \geq 0 }
\Psi_{(q , 0 , 0 , p )}^{(0, 0)}
\]
\noindent
Case $\mathrm{II}:$ For $i + j \leq n-1$ with $ i , j > 0$,
\[
\calE ^{i,j}
=
\bigoplus_{q \geq 1 , p \geq 1}
\left(
\Psi_{(q , j , i , p )}^{(i, j)}
\oplus
\Psi_{(q , j , i-1 , p )}^{(i, j)}
\oplus
\Psi_{(q , j-1 , i , p )}^{(i, j)}
\oplus
\Psi_{(q , j-1 , i-1 , p )}^{(i, j)}
\right).
\]
\noindent
Case $\mathrm{III}:$ For $1 \leq i \leq n-1$,
\[
\calE ^{i,0}
=
\bigoplus_{q \geq 0 , p \geq 1}
%V (q, \underline{0}, \underline{-1}_{i-k_2}, -p )
\left(
\Psi_{(q , 0 , i , p )}^{(i, 0)}
\oplus
\Psi_{(q , 0 , i-1 , p )}^{(i, 0)}
\right)
.
\]
\noindent
Case $\mathrm{IV}:$ For $1 \leq j \leq n-1$,
\[
\calE ^{0,j}
=
\bigoplus_{q \geq 1 , p \geq 0}
%V (q, \underline{1}_{j-k_1}, \underline{0}, -p )
\left(
\Psi_{(q , j , 0 , p )}^{(0, j)}
\oplus
\Psi_{(q , j-1 , 0 , p )}^{(0, j)}
\right)
.
\]
\noindent
Case $\mathrm{V}:$
For $ i + j = n $ with $i, j > 0$,
\[
\calE ^{i, j }
=
\bigoplus_{q \geq 1 , p \geq 1 }
\left(
%V (q, \underline{1}_{j -1}, \underline{-1}_{i}, -p )
\Psi_{(q , j , i-1 , p )}^{(i, j)}
\oplus
%\bigoplus_{q \geq 1 , p \geq 1 }
% V (q, \underline{1}_{j }, \underline{-1}_{i -1 }, -p )
\Psi_{(q , j-1 , i , p )}^{(i, j)}
 \oplus
%\bigoplus_{q \geq 1 , p \geq 1 }
% V (q, \underline{1}_{j -1}, 0, \underline{-1}_{i -1 }, -p )
 \Psi_{(q , j-1 , i-1 , p )}^{(i, j)}
\right)
 .
\]
\noindent
Case $\mathrm{VI}:$
\[
\calE ^{n , 0 }
=
\bigoplus_{q \geq -1 , p \geq 1 }
%V (q, \underline{-1}_{n-1}, -p )
\Psi_{(q , 0 , n-1 , p )}^{(n, 0)}
.
\]
\noindent
Case $\mathrm{VII}:$
\[
\calE ^{0, n }
=
\bigoplus_{q \geq 1 , p \geq -1 }
\Psi_{(q , n-1 , 0 , p )}^{(0, n)}.
\]

\end{proposition}

\section{The eigenvalues of the Rumin Laplacian}
\label{sec: The eigenvalues of the Rumin Laplacian}

\subsection{The action of $d_R$ and the Reeb vector field}
\label{sec:The action of the Rumin derivative and the Reeb vector field}

Setting $\omega_i :=  d z_i - z_i \partial |z|^2$ and
$\overline{\omega}_i := d \bar{z}_i - \bar{z}_i \bar{\partial} |z|^2$, we define differential forms
\begin{align*}
\alpha_{(j,0)}^{(0,0)}
:= & \, \sum_{\nu = 1}^{j +1} (-1)^{\nu -1}\bar{z}_{\nu} \overline{\omega}_1 \wedge \cdots \widehat{\overline{\omega_{\nu }}} \cdots \wedge \overline{\omega}_{j +1}
,
\\
\alpha_{(j,0)}^{(0,1)}
:= & \,
\overline{\omega}_1 \wedge \cdots \wedge \overline{\omega}_{j +1},
\\
\alpha_{(0,i)}^{(0,0)}
:= & \, \sum_{\mu = n- i+1}^{n +1} (-1)^{\mu -(n- i+1)} z_{\mu} \omega_{n - i +1} \wedge \cdots \widehat{\omega_{\mu }} \cdots \wedge \omega_{n +1}
,
\\
\alpha_{(0,i)}^{(1,0)}
:= & \, \omega_{n - i +1} \wedge\cdots \wedge \omega_{n +1}
.
\end{align*}
Following \cite{Orsted-Zhang-05},  we see that  $\Psi_{(q , j , i , p )}^{(s, t)}$ contains the following element $\psi_{(q , j , i , p )}^{(s, t)}$:
for $p , q \geq 1$, $a,b \geq 0$ and $a + b \leq 1$,
\begin{align*}
\psi_{(0 , 0 , 0 , 0 )}^{(0, 0)}
:= & \, 1
, \\ % End
\psi_{(q , j , i , p )}^{(i + a , j + b )}
:= & \, \overline{z}_1^{q-1} z_{n+1}^{p-1}
\alpha_{(0,i)}^{(a ,0)}
\wedge \alpha_{(j,0)}^{(0,b )}
/ \sqrt{2 \pi^{n+1}}
, \\ %End
%%
%%
%\psi_{(q , j , i , p )}^{(i+1, j)}
%:= & \, \overline{z_1}^{q-1} z_{n+1}^{p-1}
%\alpha_{(0,i)}^{(1,0)}
%\wedge \alpha_{(j,0)}^{(0,0)}
%, \\ % End
%%
%%
%\psi_{(q , j , i , p )}^{(i, j+1)}
%:= & \, \overline{z_1}^{q-1} z_{n+1}^{p-1}
%\alpha_{(0,i)}^{(0,0)}
% \wedge \alpha_{(j,0)}^{(0,1)}, \\ %End
%
%
\psi_{(q , j , i , p )}^{(i+1, j+1)}
:= & \,
P \tilde{\psi}_{(q , j , i , p )}^{(i+1, j+1)}
,
\end{align*}
where
$\tilde{\psi}_{(q , j , i , p )}^{(i+1, j+1)}
:=
\overline{z}_1^{q-1} z_{n+1}^{p-1}
\alpha_{(0,i)}^{(1,0)}
\wedge \alpha_{(j,0)}^{(0,1)}
/ \sqrt{2 \pi^{n+1}}
,$
\begin{align*}
\psi_{(q , j , 0 , 0 )}^{(0, j+b)}
:= & \, \overline{z}_1^{q-1}
\alpha_{(j,0)}^{(0,b)}
/ \sqrt{2 \pi^{n+1}}
,\\% End
%%
%\psi_{(q , j , 0 , 0 )}^{(0, j+1)}
%:= & \, \overline{z_1}^{q-1} \alpha_{(j,0)}^{(0,1)} , \\% End
%
\psi_{(0 , 0 , i , p )}^{(i+a, 0)}
:= & \, z_{n+1}^{p-1}
\alpha_{(0,i)}^{(a,0)}
/ \sqrt{2 \pi^{n+1}}
,\\
\psi_{(q , n -1 , 0 , -1 )}^{(0, n)}
:= & \, \overline{z}_1^{q-1}
\alpha_{(n+1,0)}^{(0,0)}
/ \sqrt{2 \pi^{n+1}}
, \\% End
\psi_{(-1 , 0 , n -1 , p )}^{(n, 0)}
:= & \, z_{n+1}^{p-1}
\alpha_{(0,n+1)}^{(0,0)}
/ \sqrt{2 \pi^{n+1}}
.%End
\end{align*}
We have used the projection $P$
in the definition of  $\psi_{(q , j , i , p )}^{(i+1, j+1)}$.
Let us calculate $P$ explicitly
(see also Remark \ref{rem:error_of_Orsted_Zhang} below).
Since
\[
2\Lambda (\omega_{\mu} \wedge \overline{\omega}_{\nu} )
 = - {\sqrt{-1}}{} z_{\mu} \overline{z}_{\nu} \quad \text{ for } \mu \not = \nu,
\]
we have
\[
2\Lambda \left(
\alpha_{(0,i)}^{(1,0)} \wedge \alpha_{(j,0)}^{(0,1)} \right)
={\sqrt{-1} (-1)^{ i +1}}
\alpha_{(0,i)}^{(0,0)} \wedge \alpha_{(j,0)}^{(0,0)}.
\]
%From \eqref{lem:estimate elements 4-2},
Thus
\[
2\Lambda^2 \tilde{\psi}_{(q , j , i , p )}^{(i+1, j+1)}
= {\sqrt{-1}} (-1)^{ i +1} \Lambda \psi_{(q , j , i , p )}^{(i, j)} = 0.
\]
By using the Lefschetz primitive decomposition,
we get
\begin{equation}
P |_{\Psi_{(q , j , i , p )}^{(i+1, j+1)} \oplus L \Psi_{(q , j , i , p )}^{(i, j)}}
=1 + \frac{1}{n - i - j } L \Lambda.
\label{lem:estimate elements 4-4}
\end{equation}

\begin{lemma}
\label{pr:relation of delQ between elements}

If ``$i+j \leq n-1$ and $p , \, q \geq 1$''
or ``$i \leq n-1$, $j=0$, $p \geq 1$ and $q = 0$'',
\begin{equation}
\begin{cases}
\partial_{\Rumin }
\psi_{(q , j , i , p )}^{(i, j)}
& =
( p+i)
\psi_{(q , j , i , p )}^{(i+1 , j)}
,
\medskip
\\
\bar{\partial}_{\Rumin }
\psi_{(p , i , j , q )}^{(j, i)}
& =
(-1)^{j}( p+i)
\psi_{(p , i , j , q )}^{(j, i+1)}
.
\end{cases}
\label{eq:relation of delQ between element 1-1}
\end{equation}
If $i+j \leq n-2$, $p ,q \geq 1$,
\begin{equation}
\begin{cases}
\hspace{1mm}
\partial_{\Rumin }
\psi_{(q , j , i , p )}^{(i, j+1)}
& = ( p+i)
\psi_{(q , j , i , p )}^{(i+1 , j+1)},
\medskip
\\
\bar{\partial}_{\Rumin }
\psi_{(p , i , j , q )}^{(j+1, i)}
& = (-1)^{j +1}( p+i)
\psi_{(p , i , j , q )}^{(j+1, i+1)}
.
\end{cases}
\label{eq:relation of delQ between element 1-2}
\end{equation}
Otherwise, $\partial_{\Rumin } \psi_{(q , j , i , p )}^{(s, t)} = 0$
and $\bar{\partial}_{\Rumin } \psi_{(p , i , j , q )}^{(s, t)} = 0$.

\end{lemma}

\begin{remark}
  Since $\Lambda \partial_b \psi_{(q , j , i , p )}^{(i, j)} = 0$ and $\Lambda \bar{\partial}_b \psi_{(q , j , i , p )}^{(i, j)} = 0$,
the operators $\partial_{\Rumin }$ and $\bar{\partial}_{\Rumin }$ in \eqref{eq:relation of delQ between element 1-1} coincide with $\partial_b$ and $\bar{\partial}_b$.
But, since $\Lambda \partial_b \psi_{(q , j , i , p )}^{(i, j+1)} \not = 0$ and $\Lambda \bar{\partial}_b \psi_{(q , j , i , p )}^{(i+1, j)} \not = 0$, this is not the case for  \eqref{eq:relation of delQ between element 1-2}.
\end{remark}

The action of  $\mathcal{L}_T$ on $\psi_{(q , j , i , p )}^{(s, t)}$ is also easy to compute.
Since
\[
2\mathcal{L}_T z_i = \sqrt{-1} z_i ,
\quad
2\mathcal{L}_T \omega_i = \sqrt{-1} \omega_i,
\]
we obtain
\begin{equation}
2\mathcal{L}_T \psi_{(q , j , i , p )}^{(s, t)}
= \sqrt{-1} (p +i -j -q) \psi_{(q , j , i , p )}^{(s, t)}
.
\label{lem:relation of T between element}
\end{equation}

\subsection{$L^2$-norms of highest weight vectors}
\label{subsec: L^2 norm of elements}

\begin{lemma}{\rm (\cite[Lemma 3.2]{Orsted-Zhang-05})}
\label{pr:estimate elements}
Let $p,q\ge1$ and set
\begin{align*}
C (q,p) &= 2^{n+1} \pi^{n+1} (q-1)!(p-1)!/(q+p+n)!,
\\
D(q)& = 2^{n+1} \pi^{n+1} (q-1)!/(q+n)!.
\end{align*}
If  $i + j \leq n-1$,
\begin{equation}
\left\|
\psi_{(q , j , i , p )}^{(i, j)}
\right\|^2
= \frac{ C (q,p) }{2^{j +i }}(q+j)(p+i).
\label{eq:estimate elements 1}
\end{equation}
If  $j > 0$ and $i + j \leq n-1$,
\begin{equation}
 \left\|
\psi_{(q , j , i , p )}^{(i+1, j)}
\right\|^2
=
 \left\|
\psi_{(p , i , j , q )}^{(j, i+1)}
\right\|^2
=
\frac{C(q,p)}{2^{j +i +1 }}(q+j)(q+ n -i)
.
\label{eq:estimate elements 2}
\end{equation}
If $i , j > 0 $ and $i + j \leq n-2$,
\begin{equation}
 \left\|
\psi_{(q , j , i , p )}^{(i+1, j+1)}
\right\|^2
=   \frac{C(q,p)}{2^{j +i +2 }} \frac{(q+n-i)(p+n-j)(n-1-i-j)}{n - i - j  }.
\label{eq:estimate elements 4}
\end{equation}
If  $0 \leq j \leq n-1$,
\begin{align}
 \left\|
\psi_{(q , j , 0 , 0 )}^{(0, j)}
\right\|^2
=
 \left\|
\psi_{(0 , 0 , j , q )}^{(j, 0)}
\right\|^2
= & \, \frac{D(q) }{2^{j } }(q+j) ,
\label{eq:estimate elements 6}
\\
 \left\|
\psi_{(q , j , 0 , 0 )}^{(0, j+1)}
\right\|^2
=
 \left\|
\psi_{(0 , 0 , j , q )}^{(j+1, 0)}
\right\|^2
= & \, \frac{D(q) }{2^{j+1 } }(n-j) .
\label{eq:estimate elements 7}
\end{align}

\end{lemma}

\begin{remark}
These formulas are different from those in  \cite{Orsted-Zhang-05} by factors in powers of $2$  due to the choice of the metric $g$.
\end{remark}

\begin{proof}
We only prove (\ref{eq:estimate elements 4}) because others were proved in Lemma 3.2 in \cite{Orsted-Zhang-05}; see also Remark \ref{rem:error_of_Orsted_Zhang}.
Since $P$ is the orthogonal projection and $\psi_{(q , j , i , p )}^{(i+1, j+1)} = P \tilde{\psi}_{(q , j , i , p )}^{(i+1, j+1)}$,
the formula \eqref{lem:estimate elements 4-4} gives
\[
\left\| \psi_{(q , j , i , p )}^{(i+1, j+1)} \right\|^2
=
\left\|
\tilde{\psi}_{(q , j , i , p )}^{(i+1, j+1)}
\right\|^2
 - \left\|
  (n - i - j )^{-1} L \Lambda \tilde{\psi}_{(q , j , i , p )}^{(i+1, j+1)}
 \right\|^2.
\]
The first term of the right-hand side can be calculated  by using the following facts: the squared norm of $\alpha_{(0,i)}^{(1,0)}$ in $g$ (see \cite[Lemma 5]{Folland-72}) is
$ \sum_{\mu = 1 }^{n -i} |z_{\mu} |^2 / 2^{i+1}$
and
\[
\int_{S} | z^{\alpha} |^2 d \vol_g
= \frac{2^{n+1} \pi^{n+1} \alpha !}{( |\alpha | +n )!}.
\]
For the second term, we can use
\[
\Lambda \tilde\psi_{(q , j , i , p )}^{(i+1, j+1)} = \frac{\sqrt{-1}}{2} (-1)^{i
+1} \psi_{(q , j , i , p )}^{(i, j)}
\]
and
\[\left\| L f \right\|^2
= (n -i-j) \left\| f \right\|^2,
\quad f \in \calE ^{i,j}
\]
to reduce it to
\[
\frac{1}{4(n - i - j) }\left\| \psi_{(q , j , i , p )}^{(i, j)}
\right\|^2.
\]
This is given by \eqref{eq:estimate elements 1}.
\end{proof}

\begin{remark}
\label{rem:error_of_Orsted_Zhang}
In \cite{Orsted-Zhang-05},
 the formula of the projection $P$, corresponding to our \eqref{lem:estimate elements 4-4}, is not corrct.
 This result in errors in the evaluation of the norm corresponding to our \eqref{eq:estimate elements 4} and the computations of the eigenvalues of the Laplacians using that formula.
 \end{remark}

\subsection{Calculation of eigenvalues}
\label{subsec: Calculation of Rumin derivative and Laplacians}
%この中に入っている既約成分をリストアップして，固有値の計算をする．
 Given $(q,j,i,p)$, we list up all $(s,t)$ such that  $\Psi_{(q , j , i , p )}^{(s, t)} \not = \{ 0 \}$
 and calculate the eigenvalues of $\Delta_\RN $ on them.
%この節では，
In this subsection,
we omit the subscripts from $\psi_{(q , j , i , p )}^{(s, t)}$, $\Psi_{(q , j , i , p )}^{(s, t)}$
and write
$\psi^{(s, t)}$, $\Psi^{(s, t)}$.

\medskip
\noindent
{\bf Case I:} $i=j=0$ and $p=q=0$

The space
is $\Psi^{(0, 0)}$,
and we have $\Delta_\RN \Psi^{(0, 0)} =0$.

\medskip
\noindent
{\bf Case II:}
 $i + j \leq n-2$, $p \geq 1$ and $q \geq 1$

The spaces
are
$\Psi^{(i, j)}$,
$\Psi^{(i+1, j)}$,
$\Psi^{(i, j+1)}$ and
$\Psi^{(i+1, j+1)}$.
Let $ \left\|
\partial_\RN
\right\|$
and
$\left\|
\bar{\partial}_\RN
\right\|$
be the norm of bounded linear operators of $\partial_\RN $ and $\bar{\partial}_\RN $.
By using Propositions \ref{pr:relation of delQ between elements} and \ref{pr:estimate elements}, we have
\begin{align*}
  \left\|
    \partial_\RN |_{
      \Psi^{(i, j)}
    }
  \right\|^2
  &=
    \frac{(p +i)^2}{n-i-j}
  \frac{
  \left\|
  \psi^{(i+1, j)}
  \right\|^2
  }{
  \left\|
  \psi^{(i, j)}
  \right\|^2
  }
   & & =
  \frac{(p + i)(q +n- i)}{2(n -i -j)}
  ,  \\
  \left\|
  \bar{\partial}_\RN
  |_{
  \Psi^{(i, j)}
  }
  \right\|^2
  & =
  \frac{(q +j)^2}{n-i-j}
  \frac{
    \left\|
      \psi^{(i, j+1)}
    \right\|^2
  }{
    \left\|
      \psi^{(i, j)}
    \right\|^2
  }
  & & =
    \frac{(q + j)(p +n- j)}{2(n -i -j)}
  , \\
  \left\|
    \partial_\RN |_{
      \Psi^{(i, j+1)}
    }
  \right\|^2
  &=
    \frac{(p +i)^2}{n-i-j-1}
    \frac{
      \left\|
      \psi^{(i+1, j+1)}
      \right\|^2
    }{
      \left\|
      \psi^{(i, j+1)}
      \right\|^2
    }
  & & =
    \frac{(p + i)(q +n -i)}{2(n -i -j)}
  ,\\
  \left\|
  \bar{\partial}_\RN |_{
  \Psi^{(i+1, j)}
  }
  \right\|^2
  &=
  \frac{(q +j)^2}{n-i-j-1}
  \frac{
    \left\|
    \psi^{(i+1, j+1)}
    \right\|^2
  }{
    \left\|
    \psi^{(i+1, j)}
    \right\|^2
  }
  & & =
    \frac{(q +j)(p +n -j)}{2(n -i -j)}
  .
\end{align*}
Therefore we can calculate the eigenvalue of $\Delta_\RN $ on $\Psi^{(i, j)} $ and $\Psi^{(i +1, j +1)} $ are
\begin{align*}
\Delta_\RN |_{\Psi^{(i, j)} }
&=
\left( \partial_\RN {}^{\adjoint } |_{
\Psi^{(i+1, j)}
}
\partial_\RN |_{
\Psi^{(i, j)}
}
+
\bar{\partial}_\RN {}^{\adjoint } |_{
\Psi^{(i, j+1)}
}
\bar{\partial}_\RN |_{
\Psi^{(i, j)}
}
\right)^2
 \\
&=
\frac{(( p +i )(q +n -i) + ( q +j )(p +n -j))^2}{4(n-i-j)^2}, \\
\Delta_\RN |_{
\Psi^{(i +1, j +1)}
}
&=
\left(
\partial_\RN |_{
\Psi^{(i, j+1)}
}
\partial_\RN {}^{\adjoint } |_{
\Psi^{(i+1, j+1)}
}
+
\bar{\partial}_\RN |_{
\Psi^{(i+1, j)}
}
\bar{\partial}_\RN {}^{\adjoint } |_{
\Psi^{(i+1, j+1)}
}
\right)^2
\\
&= \frac{(( p +i )(q +n -i) + ( q +j )(p +n -j))^2}{4(n -i -j)^2}
.
\end{align*}

We consider $(i + j +1)$-form.
Since $\Ima d_\RN $ and $\Ima d_\RN {}^{\adjoint }$ are orthogonal,
$\Psi^{(i+1, j)} \oplus \Psi^{(i, j+1)}
=
d_\RN \Psi^{(i, j)} \oplus d_\RN {}^{\adjoint } \Psi^{(i+1, j+1)} $.
Since $\Delta_\RN d_\RN =  d_\RN \Delta_\RN $ and
$\Delta_\RN d_\RN ^{\adjoint} =  d_\RN ^{\adjoint} \Delta_\RN $,
the eigenvalue of $\Delta_\RN $ on $\Psi^{(i+1, j)} \oplus \Psi^{(i, j+1)}$ is
\begin{align*}
\frac{( ( p +i )(q +n -i) + ( q +j )(p +n -j))^2}{4(n-i-j)^2}.
\end{align*}

\medskip
\noindent
{\bf Case III:}
$i \leq n-1$, $j = 0$, $p \geq 1$ and $q = 0$

The spaces
are $\Psi^{(i, 0)}$ and $\Psi^{(i+1, 0)}$.
We have
\[
\|
\partial_\RN |_{
\Psi^{(i, 0)}
}
\|^2
=
(p + i)/2.
\]
In the same way on Case II, we have the eigenvalue of $\Delta_\RN $ is
\[
(p + i)^2/4.
\]

\medskip
\noindent
{\bf Case IV:}
 $i =0$, $j \leq n-1$, $p=0$ and $q \geq 1$

The spaces
are $\Psi^{(0, j)}$ and $\Psi^{(0, j+1)}$.
Taking the conjugate of Case III,
the eigenvalue of $\Delta_\RN $ is
\[
(q + j)^2/4.
\]

\medskip
\noindent
{\bf Case V:}
 $i + j = n-1$, $p \geq 1$ and $q \geq 1$

The spaces
are
$\Psi^{(i, j)}$,
$\Psi^{(i+1, j)}$ and
$\Psi^{(i, j+1)}$.
We have
\begin{align*}
\left\|
\partial_\RN |_{
\Psi^{(i, j)}
}
\right\|^2&=
(p + i)(q +n- i)/2 , \\
\left\|
\bar{\partial}_\RN |_{
\Psi^{(i, j)}
}
\right\|^2&
=(q + j)(p +n- j)/2
.
\end{align*}
Therefore, on $\Psi^{(i, j)}$, the eigenvalue of $\Delta_\RN $ is
\[
{( ( p +i )(q +n -i) + ( q +j )(p +n -j) )^2}/{4}  .
\]

Next we consider $W = \Psi^{(i+1, j)} \oplus \Psi^{(i, j+1)}$.
We set
\[\underline{\psi}^{(s,t)} =
\psi^{(s,t)}/\| \psi^{(s,t)} \|.
\]
Let $A =\left\|
\partial_\RN |_{
\Psi^{(i, j)}
}
\right\| $ and
$B= \left\|
\bar{\partial}_\RN |_{
\Psi^{(i, j)}
}
\right\|$.
Then, we have
\[
d_\RN \underline{\psi}^{(i, j)}
=
A \underline{\psi}^{(i+1, j)}
+ B \underline{\psi}^{(i, j+1)}
\in \Ima d_\RN
,
\]
and
\begin{align*}
d_\RN d_\RN {}^{\adjoint }
(
A \underline{\psi}^{(i+1, j)}
+ B \underline{\psi}^{(i, j+1)}
)
= & \, d_\RN (A^2 +B^2 ) \underline{\psi}^{(i, j)} \\
= & \, (A^2 +B^2 )
(
A \underline{\psi}^{(i+1, j)}
+ B \underline{\psi}^{(i, j+1)}
)
.
\end{align*}
Therefore, eigenvalue of $\Delta_\RN $ on $\Ima d_\RN \Psi^{(i, j)}$ is
\[
(A^2 +B^2)^2 = {( ( p +i )(q +n -i) + ( q +j )(p +n -j) )^2}/{4}  .
\]

Let us find the eigenvalue on $d_\RN {\Psi^{(i, j)} }^{\bot }$,
which is the orthogonal complement in $W$.
We note that
\[
B \underline{\psi}^{(i+1, j)}
- A \underline{\psi}^{(i, j+1)}
 \in d_\RN {\Psi^{(i, j)} }^{\bot }.
\]

Let $C =( p +i -j -q)/2$,
$A' = C - 2 A^2 $ and $B' = C + 2 B^2$.
By (\ref{eq:D-1}) and \eqref{lem:relation of T between element},
\[
D
(
B \underline{\psi}^{(i+1, j)}
- A \underline{\psi}^{(i, j+1)}
)=  \sqrt{-1} \theta \wedge
(
A' B \underline{\psi}^{(i+1, j)}
- B' A \underline{\psi}^{(i, j+1)}
)
 .
\]
Since $D
(
A \underline{\psi}^{(i+1, j)}
+ B \underline{\psi}^{(i, j+1)}
) = 0 $,
we have
\begin{align*}
&
D {}^{\adjoint } D
(
B \underline{\psi}^{(i+1, j)}
- A \underline{\psi}^{(i, j+1)}
)
\\
= & \, \frac{(A' B)^2 + (B' A)^2 }{A^2 + B^2}
(
B \underline{\psi}^{(i+1, j)}
- A \underline{\psi}^{(i, j+1)}
)
.
\end{align*}
We note that
\begin{align*}
& \frac{(A' B)^2 + (B' A)^2 }{A^2 + B^2} \\
&=   \frac{1}{4} (q + j - i -p )^2 + ( p +i )(q +n -i)( q +j )(p +n -j).
\end{align*}
Under the condition $i+j=n-1$, it agrees with
\[
{ \left( ( p +i )(q +n -i) + ( q +j )(p +n -j) \right)^2 }/{4} .
\]
Therefore, we see that the eigenvalue on $\Ima d_\RN {\Psi^{(i, j)} }^{\bot }$ is
\[
{ \left( ( p +i )(q +n -i) + ( q +j )(p +n -j) \right)^2 }/{4} .
\]

\noindent
{\bf Case  VI:} $i= n-1$, $j =0$, $p \geq 1$ and $q = -1$

The space
is
$\Psi^{(n,0)}$.
Since there is no subspaces of $\calE ^{n-1} (S)$ corresponding to the $V ( \underline{-1}_{n}, -p )$,
we conclude $\partial_b {}^{\adjoint } \Psi^{(n,0)} = \bar{\partial}_b {}^{\adjoint } \Psi^{(n,0)} = \{ 0 \}$.
By (\ref{eq:D-1}), we have
\[
D \psi^{(n,0)}=  \theta \wedge \mathcal{L}_T \psi^{(n,0)}.
\]
Therefore, we have
\[
\Delta_\RN \psi^{(n,0)}
= (d_\RN d_\RN {}^{\adjoint } )^2 \psi^{(n,0)}
+ D^{\adjoint } D \psi^{(n,0)}
= \mathcal{L}_T {}^{\adjoint } \mathcal{L}_T \psi^{(n,0)}
,
\]
where $\mathcal{L}_T {}^{\adjoint }$ is the formal adjoint of $\mathcal{L}_T$ for the $L^2$-inner product.
By \eqref{lem:relation of T between element},
we see that the eigenvalue of $\Delta_\RN $ is
\[
{(p+n)^2}/{4}.
\]

\noindent
{\bf Case VII:} $i= 0$, $j =n-1$, $p = -1$ and $q \geq 1$

The space is $\Psi^{(0,n)}$.
Taking the conjugate of Case VI,
the eigenvalue of $\Delta_\RN $ is
\[
{(q+n)^2}/{4}.
\]

\begin{remark}
In Cases V-VII,
the eigenvalues of $D^{\adjoint} D$ were determined by  \cite{Julg-Kasparov-95}.
Their choice of highest weight vectors in $\Ker D$ and $\Ima D^{\adjoint} $ are different from ours.

\end{remark}

\section{Proof of Theorem \ref{theo:compute-precisely-kappa-on-the-sphere}}
\label{sec:The analytic torsion function}
From Theorem \ref{theo:all eigenvalues of the Rumin Laplacian on sphere}, we see that
the terms of $\kappa_\RN (s)$ in Cases II and V in Proposition \ref{theo:Irreducible-decomposition-of-Rumin-complex} cancel each other.
Thus we get
\begin{equation}
\kappa_\RN (s) = \kappa_{1} (s) + 2 \kappa_{2 } (s),
\label{eq:kappa_decomposition}
\end{equation}
where
\[
\kappa_{1} (s)
=
\sum_{k=0}^n
(-1)^{k+1} (n+1-k)
\dim \Ker \Delta_\RN
= \, -(n+1 ),
\]
which is the sum of the terms of $\kappa_\RN (s)$ in Case I,
and
\[
\kappa_{2} (s)
=
\sum_{i=0}^{n}
(-1 )^{i+1}
\sum_{p \geq 1 }
\frac{\dim V ( \underline{0}_{n-i} , \underline{-1}_{i} ,-p )}
{((p+i)/2)^{2s}}
,
\]
which is the sum of the terms of $\kappa_\RN (s)$ in Cases III and VI.
From Weyl's dimensional formula, we have
\begin{align*}
\dim V ( \underline{0}_{n-i} , \underline{-1}_{i}, -p )
= &
\frac{ p }{ p+i }
\begin{pmatrix}
n \\ i
\end{pmatrix}
\begin{pmatrix}
p+n \\ n
\end{pmatrix}
\\
= \frac{1}{n!} \sum_{l=1}^{n+1} e_{n+1-l}(&n-i , n-1-i, \ldots , -i) (p+i)^{l-1},
\end{align*}
where $e_l (X_0 , \ldots , X_n )$ are the elementary symmetric polynomials of $n+1$ variables of degree $l$.
Thus we get
\begin{align*}
\kappa_{2} (s)
= & \,
\frac{2^{2s+1} }{n!}
\sum_{i=0}^{n}
(-1)^{i+1}
\binom{n}{i} \sum_{p \geq 1 } \sum_{l=1}^{n+1}
\frac{
 e_{n+1-l} (n-i ,  \ldots , -i)
}{
(p+i)^{2s -l +1}
}
\\
= & \,
\frac{2^{2s+1}}{n!}
\sum_{i=0}^{n}
(-1)^{i+1}
\binom{n}{i}
\sum_{l=1}^{n} e_{n+1-l} (n-i , \ldots , -i)
\\
& \qquad \qquad
\cdot \left( \zeta (2s -l +1 ) - \sum_{k=1}^i k^{-(2s -l +1)} \right)
.
\end{align*}
Since  $
\sum_{l=1}^{n} e_{n+1-l} (n-i , \cdots , -i) k^l =0$ for $ k \le i,
$
the second sum in the last expression vanishes and we have
\[
\kappa_{2 } (s)
=
\frac{-2^{2s+1}}{n!}
\sum_{l=1}^{n}
c_l\,\zeta (2s -l +1 ),
\]
where
\[
c_l
 =
\sum_{i=0}^{n}
(-1)^i
\binom{n}{i}
e_{n+1-l} (n-i , \ldots , -i)
.
\]
If we set $\sigma(k)=\sum_{l=0}^{n+1}c_l\, k^l$ for
 $k\in\mathbb{N}$, then
\begin{align*}
\sigma (k)
& =
\sum_{i=0}^{n}
(-1)^i
\binom{n}{i}
\prod_{l=0}^{n}
(k +n -i -l)
\\
& =
\sum_{i=0}^{n}
(-1)^i
\binom{n}{i}
\frac{d^{n+1}t^{k+n-i}}{dt^{n+1}} \Big |_{t=1}
=
\frac{d^{n+1}t^k(t-1)^n}{dt^{n+1}} \Big |_{t=1}
.
\end{align*}
It follows that  $\sigma (k)=(n+1)! \, k$ and hence
$c_l=0$ except $c_1=(n+1)!$.
Summing up, we conclude
\begin{align*}
\kappa_\RN (s)
=
- (n+1) (1+ 2^{2s+1} \zeta (2s ))
.
\end{align*}
Using $\zeta (0)= - 1/2$ and  $\zeta ' (0) = - (\log 2 \pi ) /2$,
we get
\begin{align*}
\kappa_\RN (0)
&=
- (n+1)
(
1
+2 \zeta (0 )
)
= \,
0
,
\\
\kappa_\RN ' (0)
&=  \, 2 (n+1)  \log 4 \pi
\end{align*}
as claimed.

\end{document}